\documentclass[11pt,oneside]{article}
\usepackage[utf8]{inputenc}
\usepackage{amsmath} 
\usepackage{amsthm}
\usepackage{textcomp}
\usepackage{amssymb}
\usepackage[bottom=1cm]{geometry}

\usepackage{listings}  
\usepackage{tikz}
\usepackage{mathdots}
\usetikzlibrary{matrix}
\usepackage{vmargin}
\usepackage{thmtools}

\usepackage{mathtools}
\usepackage{dsfont}
\usepackage{multicol}
\usepackage[english]{babel}

\usepackage{caption}

\theoremstyle{plain}
\usepackage{graphicx}
\usepackage{subcaption}
\theoremstyle{plain}
\newtheorem{teo}{}[section] 

\newtheorem{theorem}[teo]{Theorem}

\newtheorem{defi}[teo]{Definition}
\newtheorem{remark}[teo]{Remark}
\newtheorem{example}[teo]{Example}

\newtheorem{prop}[teo]{Proposition}
\newtheorem{lem}[teo]{Lemma}

\newcommand\blfootnote[1]{%
  \begingroup
  \renewcommand\thefootnote{}\footnote{#1}%
  \addtocounter{footnote}{-1}%
  \endgroup
}

\title{Topological realizations of groups in Alexandroff spaces}
\author{Pedro J. Chocano, Manuel A. Mor\'on and Francisco R. Ruiz del Portal}
\date{}
\begin{document}

\maketitle

\begin{abstract}
We prove that every group can be realized as the homeomorphism group and as the group of (pointed) homotopy classes of (pointed) self-homotopy equivalences of infinitely many non-homotopy-equivalent Alexandroff spaces.
\end{abstract}

\section{Introduction}
\blfootnote{2010  Mathematics  Subject  Classification: 55P10, 55P99, 06A06.}
\blfootnote{Keywords: Automorphisms, homotopy equivalence, Alexandroff spaces, posets.}
\blfootnote{This research is partially supported by Grants MTM2015-63612-P, PGC2018-098321-B-100 and BES-2016-076669 from Ministerio de Ciencia, Innovación y Universidades (Spain).}

Growing out of foundational papers on the subject, such as \cite{mccord1966singular} and \cite{stong1966finite}, the theory of finite topological spaces has seen a substantial development in the last two decades. This is mainly due to the interest of researchers in applying computational methods to, for instance, dynamical systems \cite{barmak2020Lefschetz,mrozek2017conley} or approximation of spaces \cite{moron2008connectedness,mondejar2020reconstruction}. Two main monographs on the issue of the algebraic aspects of the topology of finite spaces are \cite{may1966finite} and \cite{barmak2011algebraic}. The first one, due to J.P. May, is the result of some REU programs developed by the author at the University of Chicago. The other one is, essentially, the Ph. D. thesis of J.A. Barmak (under the supervision of E.G. Minian).

Finite topological spaces are a particular case of more general topological spaces that were introduced by P.S. Alexandroff  \cite{alexandroff1937diskrete}. An Alexandroff space is a topological space with the property that the arbitrary intersection of open sets is open. Some of the results given by M.C. McCord  \cite{mccord1966singular} can be rephrased in the following way:
\textit{Given a polyhedron $X$, there exists an Alexandroff space $\mathcal{X}(X)$ such that the homotopy groups of $X$ and $\mathcal{X}(X)$ are the same.}
Then, it can be deduced that every group can be realized as the fundamental group of an Alexandroff space. In addition, a similar statement can be obtained for abelian groups and higher homotopy groups. 

Alexandroff spaces can also be used to realize finite groups as homeomorphism groups, e.g.,  \cite{barmak2009automorphism,birkhoff1936order,thornton1972spaces}. The problem of realizing a group as the group of homeomorphisms of a topological space, i.e., the realizability problem for the topological category ($Top$), has been widely studied. We have only cited some references related to finite topological spaces. In particular, in \cite{barmak2009automorphism}, J.A. Barmak and E.G. Minian focused on giving, for any finite group $G$, a finite topological space with the lowest possible cardinality having $G$ as the corresponding group of homeomorphisms. 

Given a topological space $X$, there are more groups that can be related to it. For instance, the group of (pointed) homotopy classes of (pointed) self-homotopy equivalences, that is to say, the group of automorphisms of $X$ $((X,x)$ where $x\in X$ $)$ as an object of $HTop$ ($HTop_*$). $HTop$ denotes the homotopy category of topological spaces and $HTop_*$ denotes the pointed homotopy category of pointed topological spaces.  In $HPol$ ($HPol_*$), the full subcategory of $HTop$ ($HTop_*$) whose objects are all (pointed) topological spaces having the homotopy type of a (pointed) polyhedron and whose morphisms are the (pointed) homotopy classes of (pointed) continuous maps, the corresponding realizability problem is of interest in the literature. This problem has appeared in many papers for over fifty years, for example  \cite{arkowitz1990group,felix2010problems,kahn1976realization,rutter2006spaces}, and has been placed as the first problem to solve in \cite{arkowitz2001problems}, a list of open problems about groups of self-homotopy equivalences. In this direction, a complete answer for finite groups and pointed topological spaces  was obtained by C. Costoya and A. Viruel \cite{costoya2014every}.
\begin{theorem}[\cite{costoya2014every}]\label{costoyaViruel}
Every finite group $G$ can be realized as the group of self-homotopy equivalences of infinitely many (non-homotopy-equivalent) rational elliptic spaces $X$.
\end{theorem}
In a recent paper \cite{costoya2014notes}, the free case has been completely solved using tools of highly algebraic character and for Eilenberg-MacLane spaces.



In conclusion, we want to show that Alexandroff spaces are a good tool to solve realization problems in topological categories such as $Top$, $HTop$ and $HTop_*$. Specifically, we are interested in giving a positive answer to the following question, where $\mathcal{C}$ is one of the previous categories. \textit{Given an arbitrary group $G$, is there an object $X$ in $\mathcal{C}$ such that its group of automorphisms is isomorphic to $G$?}

We state the three main results of the paper. The idea for the first one is to generalize the construction for finite groups of J.A. Barmak and E.G. Minian made in \cite{barmak2009automorphism}.

\begin{theorem}\label{thm:Aut(X_G)isoG}
Every group can be realized as the group of homeomorphisms of an Alexandroff space.
\end{theorem}

In \cite{stong1966finite}, R.E. Stong introduced the concepts of core and minimal finite space in order to study the homotopy type of finite topological spaces. The concept of minimal space plays a central role herein because it has the property that a self-homotopy equivalence is indeed a homeomorphism. In \cite{kukiela2010homotopy}, M.J. Kukie{\l}a extended some of these concepts and results to Alexandroff spaces, which is our framework. Then, we modify the space constructed in the proof of Theorem \ref{thm:Aut(X_G)isoG} taking into account the previous concepts so as to obtain the following result.
\begin{theorem}\label{thm:equivInfinite}
Every group can be realized as the group of self-homotopy equivalences of infinitely many (non-homotopy-equivalent) Alexandroff spaces.
\end{theorem}

The starting point to get the third result is the construction made in the proof of Theorem \ref{thm:equivInfinite}. 





\begin{theorem}\label{cor:pointedHomotopy} Every group can be realized as the group of pointed homotopy classes of pointed self-homotopy equivalences of infinitely many (non-homotopy-equivalent) Alexandroff spaces.
 \end{theorem}

The organization of the paper is as follows. In Section \ref{section:preliminares}, we introduce the basic definitions and theorems from the literature that we will use in subsequent sections. In section \ref{section:constructionX_G}, given a group $G$, we provide a method to obtain three Alexandroff spaces ($X_G$, $\overline{X}_G$ and $\overline{X}_G^*$) that are the candidates to solve the problem of realizability for the categories $Top$, $HTop$ and $HTop_*$. Moreover, we present one example in detail to illustrate each construction. In Section \ref{section:Automorphism(X_G)}, we show that the group of homeomorphisms of $X_G$ is isomorphic to $G$, solving the problem of realizability for the topological category. 
In section \ref{section:selfhomotopyequivalences}, we prove that every group can be realized as the group of (pointed) homotopy classes of (pointed) self-homotopy equivalences of infinitely many non-homotopy-equivalent Alexandroff spaces, i.e., Theorem \ref{thm:equivInfinite} and Theorem \ref{cor:pointedHomotopy}. The proof of Theorem \ref{thm:equivInfinite} is divided into three auxiliary lemmas and uses the main result of Section \ref{section:Automorphism(X_G)}. The proof of Theorem \ref{cor:pointedHomotopy} is similar to the proof of Theorem \ref{thm:equivInfinite}. In section \ref{Section:McCord}, we study some properties of the space $\overline{X}_G$ and its McCord complex $\mathcal{K}(\overline{X}_G)$. 


\section{Preliminaries}\label{section:preliminares}

We introduce concepts that will be used in subsequent sections. 
\begin{defi} An Alexandroff space $X$ is a topological space with the property that the arbitrary intersection of open sets is open. Furthermore, if $x\in X$, $U_x$ denotes the open set that is defined as the intersection of all open sets containing $x$. Analogously, $F_x$ denotes the closed set that is given by the intersection of all closed sets containing $x$.
\end{defi}
\begin{defi} Given a partially ordered set (poset) $(X,\leq)$, a lower (upper) set $S$ is a subset of $X$ such that if $x\in S$ and $y\leq x$ ($x\leq y$), then $y\in S $.
\end{defi}
\begin{theorem}[\cite{alexandroff1937diskrete}]\label{thm:Alex}
For a poset $(X,\leq)$, the family of lower (upper) sets of $\leq$ is a $T_0$ topology on $X$, that makes $X$ an Alexandroff space. For a $T_0$ Alexandroff space, the relation $x\leq_{\tau} y$ if and only if $U_x\subset U_y$ ($U_y \subset U_x$) is a partial order on $X$. Moreover, the two associations relating $T_0$ topologies and partial orders are mutually inverse. 
\end{theorem} 
\begin{remark} $U_x$ can also be seen as the set $\{y\in X|y\leq x\}$. Similarly, $F_x$ can also be seen as the set $\{y\in X|y\geq x\}$. Furthermore, $U_x$ and $F_x$ are contractible spaces. 
\end{remark}

The second order and topology defined above, the ones in parenthesis, are usually called the opposite order and opposite topology. From now on, we treat posets and Alexandroff spaces as the same objects without explicit mention. We assume that the orders and topologies are not the opposite ones. We also assume that all the Alexandroff spaces that appear throughout this section satisfy the separation axiom $T_0$.

By Theorem \ref{thm:Alex}, some topological notions can be expressed  using partial orders. For instance, let $X$ and $Y$ be  two Alexandroff spaces, $f:X\rightarrow Y$ is a continuous function if and only if $f$ is order-preserving. If $X$ and $Y$ are finite topological spaces, $f,g:X\rightarrow Y$ are homotopic if and only if there exists a sequence $f_0,...,f_n$ of continuous maps from $X$ to $Y$ such that $f(x)=f_0(x)\leq f_1(x)\geq f_2(x)\leq...f_n(x)=g(x)$ for every $x$ in $X$. In the case that $X$ and $Y$ are Alexandroff spaces and $f(x)\geq g(x)$ for every $x\in X$,  $f$ is homotopic to $g$. Moreover, if $x_0,x_{1},...,x_{n}$ are points in an Alexandroff space $X$ such that $x_i $ is comparable to $x_{i+1}$ for every $i=0,...,n-1$, there exists a path from $x_0$ to $x_n$. See \cite{barmak2011algebraic,may1966finite} for more details.

Given a poset $(X,\leq)$, we denote $x\prec y$ ($x\succ y $) if and only if $x<y$ ($x>y$) and there is no $z$ such that $x<z<y$ ($x>z>y$). If $x\prec y$ ($x\succ y $), we say that $x$ is an immediate predecessor of $y$ ($x$ is an immediate successor of $y$). We also denote by $maximal (X)$  ($minimal(X)$) the set of maximal (minimal) elements of $X$, while $max(X)$ denotes the maximum of $X$. Inequalities of the form $y\leq x<max(X)$, where $y,x\in X$, will appear in subsequent sections. If $max(X)$ does not exist, we assume that $y\leq x<max(X)$ means $y\leq x$.
\begin{prop}\label{prop:sucesores} Let $f:X\rightarrow Y$ be a homeomorphism between two Alexandroff spaces. If $x,y\in X$ are such that $x\prec y$, then $f(x)\prec f(y)$.
\begin{proof}
We argue by contradiction, suppose that there exists $z\in Y$ with $f(x)<z<f(y)$. Applying $f^{-1}$, we get $x<f^{-1}(z)<y$, which leads to a contradiction.
\end{proof}
\end{prop}

R.E. Stong \cite{stong1966finite} provided a method to classify finite spaces by their homotopy type. Despite the fact that the following definition was enunciated for finite spaces, we enunciate it for Alexandroff spaces.
\begin{defi} Let $X$ be an Alexandroff space.
 \begin{itemize}
\item $x\in X$ is linear, or up beat point following modern terminology, if there exists $y>x$ with the property that for every $z>x$ we have $z\geq y$. Equivalently, $F_x\setminus\{ x\}$ has a minimum.
\item $x\in X$ is colinear, or down beat point following modern terminology, if there exists $y<x$ with the property that for every $z<x$ we have $y \geq z$. Equivalently, $U_x\setminus\{ x\}$ has a maximum.
\end{itemize}
Moreover, we say that a finite topological space $X$ is a minimal finite space if $X$ is $T_0$ and has no linear or colinear points (beat points). A core of a finite space $X$ is a strong deformation retract of it which is a minimal finite space.
\end{defi}
If we remove a colinear or a linear point (beat point) of a finite topological space $X$, the homotopy type of $X$ does not change \cite{stong1966finite}. Therefore, we can remove beat points one by one until we get a minimal finite space, which is also a strong deformation retract of $X$.  The following result shows that the core obtained is unique up to homeomorphism.

\begin{theorem}[\cite{stong1966finite}]\label{thm:stongClassification}
Let $X$ and $Y$ be minimal finite spaces. Then, $X$ is homotopy equivalent to $Y$ if and only if $X$ is homeomorphic to $Y$.
\end{theorem}

The notion of core introduced by R.E. Stong and some results for finite spaces were generalised to Alexandroff spaces by M.J. Kukie{\l}a \cite{kukiela2010homotopy}. 
We recall two definitions.
\begin{defi} Let $X$ be an Alexandroff space, $r:X\rightarrow X$ is a comparative retraction if $r$ is a retraction in the usual sense and for all $x\in X$, $r(x)\leq x$ or $r(x)\geq x$. The class of all comparative retractions is denoted by $\mathcal{C}$. The space $X$ is called a $\mathcal{C}$-core if there is no other retraction $r:X\rightarrow X $ in $\mathcal{C}$ other than the identity $id_X$.
\end{defi}
\begin{remark}\label{rem:beatEntoncesNoCore} If $X$ is an Alexandroff space with a beat point $x$, then $X$ is not a $\mathcal{C}$-core. If we assume that $x$ is an up beat point, $F_x\setminus \{x\}$ has a minimum $x'$. We can define $r:X\rightarrow X$ given by $r(x)=x'$ and $r(y)=y$ for every $y\in X\setminus \{ x\}$. It is clear that $r$ is a comparative retraction and it is not the identity. A similar argument can be used for a down beat point.
\end{remark}
\begin{defi}\label{def:locallycore} We say a $\mathcal{C}$-core $X$ is locally a core if for every $x\in X$ there exists a finite set $A_x\subset X$ containing $x$ such that for every $y \in A_x$, then  $|A_x\cap maximal(\{ z\in X |z<y \})|\geq 2$ if $y$ is not minimal in $X$ and $|A_x\cap minimal(\{ z\in X |z>y \})|\geq 2$ if $y$ is not maximal in $X$.
\end{defi}
Given two topological spaces $X$ and $Y$, the space of continuous maps from $X$ to $Y$ equipped with the compact-open topology will be denoted by $C(X,Y)$.
\begin{theorem}[\cite{kukiela2010homotopy}]\label{thm:kukielaequivIdentity}
If $X$ is locally a core, then there is no map in $C(X,X)$ homotopic to the identity $id_X$ other than $id_X$.
\end{theorem}

Given a topological space $X$, $Aut(X)$ denotes the group of homeomorphisms of $X$ and $\mathcal{E}(X)$ denotes the group of homotopy classes of self-homotopy equivalences of $X$.

\begin{lem}\label{lem:clave} If $X$ is locally a core, then $Aut(X)$ and $\mathcal{E}(X)$ are isomorphic.
\begin{proof}
We define $\varphi:Aut(X)\rightarrow \mathcal{E}(X)$ given by $\varphi(f)=[f]$, where $[f]$ denotes the homotopy class of $f$. $\varphi$ is clearly well-defined and a homomorphism of groups. If $\varphi(f)$ is the homotopy class of the identity map, we get that $f$ is the identity map by Theorem  \ref{thm:kukielaequivIdentity}. Therefore, $\varphi$ is a monomorphism of groups. If $f$ is a self-homotopy equivalence of $X$, there exists a continuous map $g:X \rightarrow X$ such that $f\circ g$ and $g\circ f$ are homotopic to the identity map $id_X$. By  Theorem \ref{thm:kukielaequivIdentity}, we obtain that $g\circ f=id_X$ and $f\circ g=id_X$. Thus, it can be deduced that $f$ is indeed a homeomorphism. From here, we can obtain that $\varphi$ is an isomorphism of groups.
\end{proof}
\end{lem}
\begin{remark}\label{rem:AutEqSonIguales} If $X$ is locally a core, it can be deduced easily that each homotopy class of the group $\mathcal{E}(X)$ contains exactly one element. We will refer and treat an element $[f]$ of $\mathcal{E}(X)$ just as $f\in Aut(X)$ and vice versa. Therefore, we can identify both groups, so $Aut(X)=\mathcal{E}(X)$.
\end{remark}

\begin{defi} A weak homotopy equivalence is a map between topological spaces  which induces isomorphisms on all homotopy groups. Furthermore, it is said that two topological spaces $X,Y$ are weak homotopy equivalent (or they have the same weak homotopy type) if there exists a sequence of spaces $X=X_0,X_1,...,X_n=Y$ such that there exist weak homotopy equivalences $X_i\rightarrow X_{i+1}$ or $X_{i+1}\rightarrow X_i$ for every $0\leq i\leq n-1$. 
\end{defi}

In \cite{mccord1966singular}, M.C. McCord studied the weak homotopy type of Alexandroff spaces using simplicial complexes. The key to get the most important result in that paper relies in the following theorem, which is somehow an adaptation of a theorem by A. Dold and R. Thom \cite{dold1958quasifaserungen}.


\begin{defi} An open cover $\ \mathcal{U}$ of a space $B$ will be called basis-like if whenever $U,V\in \mathcal{U}$ and $x\in U\cap V$, there exists $W\in \mathcal{U}$ such that $x\in W\subset U\cap V$.
\end{defi}
Every Alexandroff space $X$ admits a basis-like open cover given by $\mathcal{U}=\{ U_x|x\in X\}$. 
\begin{theorem}[\cite{mccord1966singular}]\label{thm:teorema6McCord} Suppose $p$ is a map of a space $E$ into a space $B$ for which there exists a basis-like open cover $\mathcal{U}$ of $B$ satisfying the following condition: for each $U\in \mathcal{U}$, the restriction $p_{|p^{-1}(U)}:p^{-1}(U)\rightarrow U$ is a weak homotopy equivalence. Then $p$ itself is a weak homotopy equivalence.
\end{theorem}

\begin{defi} Let $X$ be an Alexandroff space. It can be considered the order complex (or McCord complex) $\mathcal{K}(X)$, i.e., the vertices of the complex are the points of $X$ and the simplices are the finite, totally ordered subsets of $X$. The geometric realization of $\mathcal{K}(X)$ is denoted by $|\mathcal{K}(X)|$.
\end{defi}

Let $X$ be an Alexandroff space, for every $u\in |\mathcal{K}(X)| $ we have that $u$ is contained in a unique open simplex $(x_0,...,x_r)$, where $x_0<...<x_r$. Then, $f_X:|\mathcal{K}(X)|\rightarrow X$ is defined by $f_X(u)=x_0$. M.C. McCord showed that $f_X$ is continuous and has the property that $|\mathcal{K}(U_x)|$ is a deformation retract of $f^{-1}_X(U_x)$ and contractible. From here, 
\begin{theorem}[\cite{mccord1966singular}]\label{thm:McCord}
There exists a correspondence that assigns to each Alexandroff space $X$ a simplicial complex $\mathcal{K}(X)$ and a weak homotopy equivalence $f_X:|\mathcal{K}(X)|\rightarrow X$.
\end{theorem}

Hasse diagrams are usually defined for finite posets, but they can also be defined for some non-finite posets. Let $X$ be an Alexandroff space, it will be said that $X$ is locally finite if for every $x\in X$ there are finitely many points smaller than $x$. If $X$ is a locally finite Alexandroff space, the Hasse diagram of $X$ is a directed graph $H(X)$, where the vertices are the points of $X$ and there is an edge between two vertices $x$ and $y$ if and only if $x\prec y$ (or $x\succ y$). The direction of the edge goes from the lower element to the upper element. In the subsequent Hasse diagrams, we will omit the orientation and we will assume an upward orientation. On the other hand, not every Alexandroff space can be represented by a Hasse diagram, for instance, the real numbers with the usual order. An example of a Hasse diagram for a non-finite poset can be found in Example \ref{ex:infiniteposet}.

\section{Construction of $X_G$, $\overline{X}_G$ and $\overline{X}^*_G$ }\label{section:constructionX_G}

Firstly, given a group $G$, we construct an Alexandroff space $X_G$ that will be used to prove Theorem \ref{thm:Aut(X_G)isoG}. The construction is analogous to the one given by J.A. Barmak and E.G. Minian for finite groups \cite{barmak2009automorphism}. 

Let $G$ be a group and $S'$ be a set of non-trivial generators of $G$, i.e., the identity element is not in $S'$. By the well-ordering principle, we can take a well-order on the set $S'$. If $S'$ is not finite and contains a maximum such that $max(S')$ does not have an immediate predecessor, we extend the well-order defined on $S'\setminus\{ max(S')\}$ as a subset of $S'$ to $S'$  declaring that $max(S')<\alpha$ for every $\alpha\in S'\setminus\{max(S')\}$. It is easy to check that this order on $S'$ is indeed a well-order. Now, $S'$ does not have a maximum, otherwise, we get a contradiction with the fact that $S'$ with the previous order satisfies that $max(S')$ does not have an immediate predecessor. 

Thus, $S'$ with the well-order considered before satisfies only one of the following three properties:
\begin{itemize}
\item $S'$ is finite.
\item $S'$ is not finite and does not have a maximum.
\item $S'$ is not finite and has a maximum such that $max(S')$ has an immediate predecessor.
\end{itemize}

We consider $S=S'\cup \{-1,0 \}$, where we are assuming that $-1,0\notin S $. We extend the well-order of $S'$ to $S$, for every $\alpha\in S'$ we declare $-1<0<\alpha$. Finally, we consider $X_G=G \times S $ with the following relations:
\begin{itemize}
\item $(g,\beta)<(g,\gamma)$ if $-1\leq \beta<\gamma$ where $g\in G$ and $\beta,\gamma\in S$.
\item $(g\beta,-1)<(g,\gamma)$ if $0<\beta\leq\gamma$ where $g\in G$ and $\beta,\gamma\in S$.
\end{itemize}

It is trivial to check that $X_G$ with the previous relations is a partially ordered set. Therefore, by Theorem \ref{thm:Alex}, it is a $T_0$ Alexandroff space.
\begin{figure}[h]
  \centering
    \includegraphics[scale=1]{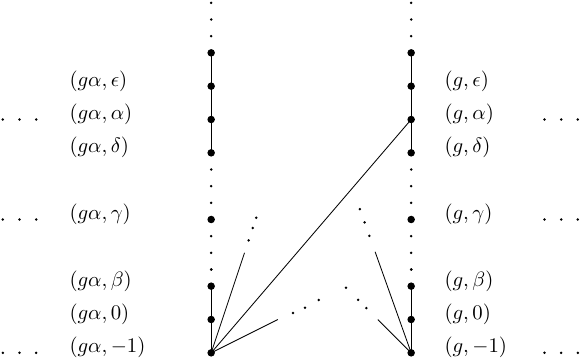}
     \caption{Schematic Hasse diagram of $X_G$.}
     \label{figure:schemeHasseDiagram}         
\end{figure}



If $X_G$ is locally a core, we get that $Aut(X_G)$ is isomorphic to $\mathcal{E}(X_G)$ by Lemma \ref{lem:clave}. Therefore, it is enough to show that $Aut(X_G)$ is isomorphic to $G$ so as to obtain that every group can be realized as the group of homotopy classes of self-homotopy equivalences of an Alexandroff space. But,
$X_G$ is far from being locally a core because it contains beat points, Remark \ref{rem:beatEntoncesNoCore}. Concretely, every point of the form $(g,\beta)$, where $g\in G$ and $0\leq \beta< max(S)$, is a beat point. $F_{(g,\beta)}\setminus\{ (g,\beta)\}$ can be seen as a subset of $S$ if we ignore the first coordinate ($F_{(g,\beta)}\setminus\{ (g,\beta)\}=\{(g,\alpha)|\alpha>\beta\}$). Then, there exists a minimum in $F_{(g,\beta)}\setminus\{ (g,\beta)\}$, which means that $(g,\beta)$ is an up beat point. We need to add points to $X_G$ so as to get a good candidate $\overline{X}_G$ to be locally a core. At the same time, we also want to modify $X_G$ in order to get that $Aut(\overline{X}_G)$ is isomorphic to $Aut(X_G)$. By Lemma \ref{lem:clave} and showing that $Aut(X_G)$ is isomorphic to $G$, we would get the desired result.

The idea of the construction is to change the state of the beat points $(g,\beta)$ of $X_G$. If we add the point $B_{(g,\beta)}$ to $X_G$ with the relation $B_{(g,\beta)}>(g,\beta)$, we get easily that $(g,\beta)$ is not a beat point. On the other hand, the new point added is now a beat point. To solve this situation, we can add the point $C_{(g,\beta)}$ with the relation $C_{(g,\beta)}<B_{(g,\beta)}$. We have that $B_{(g,\beta)}$ is not a beat point but $C_{(g,\beta)}$ is a beat point. We argue following the previous method so as to get $S_{(g,\beta)}$, i.e., we add $A_{(g,\beta)}$ and $D_{(g,\beta)}$ with $C_{(g,\beta)}<A_{(g,\beta)}>D_{(g,\beta)}<(g,\beta)$. Now, the new points added and $(g,\beta)$ are not beat points.

\begin{figure}[h]
  \centering
 \includegraphics[scale=0.8]{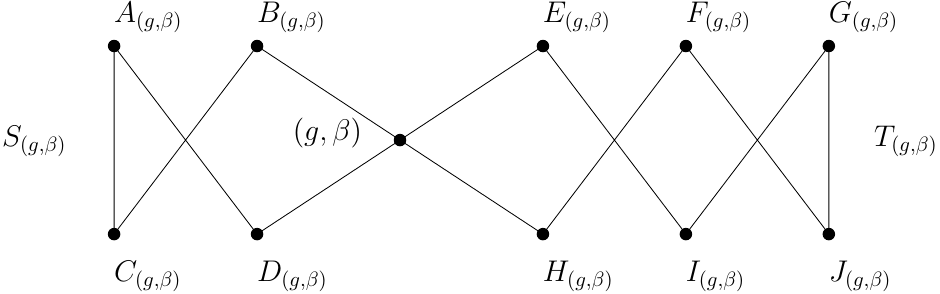}
     \caption{Hasse diagram of $S_{(g,\beta)}\cup T_{(g,\beta)}$.}
     \label{figure:hassediagramAspas}         
\end{figure}

We construct $T_{(g,\beta)}$ following the idea of not introducing beat points, Figure \ref{figure:hassediagramAspas}. It is also important to get that $T_{(g,\beta)}$ is not homeomorphic to $S_{(g,\beta)}$. If $S_{(g,\beta)}$ is homeomorphic to $T_{(g,\beta)}$, we cannot expect to obtain that $Aut(X_G)$ is isomorphic to $Aut(\overline{X}_G)$ because we introduce more homeomorphisms, particularly, the homeomorphism $f$ that satisfies $f(T_{(g,\beta)})=S_{(g,\beta)}$, $f(S_{(g,\beta)})=T_{(g,\beta)}$ and keeps the rest of points fixed. In addition, the group of homeomorphisms of $S_{(g,\beta )}\cup T_{(g,\beta)}$ is trivial as we will see in Lemma \ref{lem:SyTsonasimetricos}. Now, it is easy to find the set $A_x$ asked in Definition \ref{def:locallycore} for $x=(g,\beta)$, where $0\leq \beta<max(S)$, $A_x=S_{(g,\beta)}\cup T_{(g,\beta)}$. It is also easy to find $A_x$ for the rest of the points as we will see in the proof of Lemma \ref{lem:E(Xbarra_G)Aut(Xbarra_G)}. However, if we do not consider the space $T_{(g,\beta)}$ in the construction of $\overline{X}_G$, the existence of a set $A_x$, for every point $x$ in $\overline{X}_G$, satisfying the conditions asked in Definition \ref{def:locallycore} is not guaranteed. Suppose $G$ is a group and $S'$ a non-finite countable generating set of $G$. We can denote the elements of $S'$ by $s_i$ with $i\in \mathbb{N}$. We can define a well-order on $S'$ given by $s_i<s_j$ if and only if $i<j$. If we consider $X=X_G\cup (\bigcup_{\alpha\in S\setminus\{-1 \}} S_{(g,\alpha)})$ and the point $x=(g,0)$ for some $g\in G$, a set $A_x$ satisfying the conditions of Definition \ref{def:locallycore} does not exist. If it exists, $(g,s_i)\in A_x$ for every $i\in \mathbb{N}$, which means that 
$A_x$ cannot be a finite set. We prove the last assertion. By construction, for every $j\in \mathbb{N}$, $(g,s_j)$ is not a maximal point. Furthermore, $minimal(\{ z\in X|z>(g,0) \})=\{(g,s_1),B_{(g,0)} \}$. Then, $(g,s_1),B_{(g,0)}\in A_x$ because it is needed that $|A_x\cap minimal(\{ z\in X|z>(g,0) \})|\geq 2 $. Moreover, $minimal(\{ z\in X|z>(g,s_1) \})=\{(g,s_2),B_{(g,s_1)} \}$, so $(g,s_2),B_{(g,s_1)}\in A_x$. In general, it can be argued inductively that $(g,s_{n+1}),B_{(g,s_n)}\in A_x$ for every $n\in \mathbb{N}$ because $minimal(\{ z\in X|z>(g,s_n) \})=\{(g,s_{n+1}),B_{(g,s_n)} \}$.
 
If $(g,\beta)\in X_G$, where $0\leq \beta<max(S)$ and $g\in G$, we consider $S_{(g,\beta)}$ and $T_{(g,\beta)}$ in the following way: $S_{(g,\beta)}=$ $\{A_{(g,\beta)},$ $B_{(g,\beta)},$ $C_{(g,\beta)},D_{(g,\beta)},(g,\beta)\}$ and  $T_{(g,\beta)}=\{ E_{(g,\beta)},$  $F_{(g,\beta)},$ $G_{(g,\beta)},$ $H_{(g,\beta)}, I_{(g,\beta)},$ $J_{(g,\beta)},$ $(g,\beta) \}$. Finally, we consider $$\overline{X}_G= X_G\cup(\bigcup_{\substack{{(g,\beta)\in G\times S} \\ 0\leq \beta <max(S)}}(S_{(g,\beta)}\cup T_{(g,\beta)}))$$
with the following relations:
\begin{enumerate}
\item $(g,\beta)<(g,\gamma)$ if $-1\leq \beta<\gamma$, where $g\in G$ and $\beta,\gamma\in S$.
\item $(g\beta,-1)<(g,\gamma)$ if $0<\beta\leq\gamma$, where $g\in G$ and $\beta,\gamma\in S$.
\item $(g,\beta)>D_{(g,\gamma)},H_{(g,\gamma)}$ if $-1<\gamma< \beta$, where $g\in G$ and $\beta,\gamma\in S$.
\item $(g,\beta)>D_{(g,\beta)},H_{(g,\beta)}$ if $0\leq \beta <max(S)$, where $g\in G$ and $\beta\in S$.
\item $(g,\beta)<E_{(g,\gamma)},B_{(g,\gamma)}$ if $-1\leq \beta<\gamma<max(S)$, where $g\in G$ and $\beta,\gamma\in S$.
\item $(g,\beta)<E_{(g,\beta)},B_{(g,\beta)}$ if  $0\leq \beta <max(S)$, where $g\in G$ and $\beta\in S$.
\item  $A_{(g,\beta)}>C_{(g,\beta)},D_{(g,\beta)}$ and $B_{(g,\beta)}>C_{(g,\beta)}$ if $0\leq \beta <max(S)$, where $g\in G$ and $\beta\in S$.
\item $B_{(g,\beta)}>D_{(g,\gamma)},H_{(g,\gamma)}$ if  $-1<\gamma\leq \beta< max(S)$, where $g\in G$ and $\beta,\gamma\in S$.
\item $E_{(g,\beta)}>D_{(g,\gamma)},H_{(g,\gamma)}$ if  $-1<\gamma\leq \beta <max(S)$, where $g\in G$ and $\beta,\gamma\in S$.
\item $ G_{(g,\beta)}>I_{(g,\beta)},J_{(g,\beta)}$ if $0\leq \beta< max(S)$, where $g\in G$ and $\beta\in S$.
\item $F_{(g,\beta)}>H_{(g,\beta)},J_{(g,\beta)}$ and $E_{(g,\beta)}>I_{(g,\beta)}$ if $0\leq \beta <max(S)$, where $g\in G$ and $\beta\in S$.

%

\end{enumerate}

It is not complicated to check that $\overline{X}_G$ with the previous relations is a partially ordered set, which means that $\overline{X}_G$ is a $T_0$ Alexandroff space.

To conclude, we only need to add one extra point $*$ to $\overline{X}_G$ so as to obtain the pointed case for Theorem \ref{thm:equivInfinite}, i.e., Theorem \ref{cor:pointedHomotopy}. The point $*$ will play the role of a fixed point for every self-homotopy equivalence. We consider $\overline{X}^*_G$ as the union of $\overline{X}_G$ and $*$. We extend the partial order of $\overline{X}_G$ to $\overline{X}_G^*$ declaring that $*>(g,-1)$ for every $g\in G$.

\begin{example} Let us consider the Klein-four group $\mathbb{Z}_2\oplus \mathbb{Z}_2$, we take $S'=\{a,b\}$ as a set of non-trivial generators, where $a=(1,0)$, $b=(0,1)$. We declare $a<b$. The Hasse diagrams of $X_{\mathbb{Z}_2\oplus \mathbb{Z}_2}$, $\overline{X}_{\mathbb{Z}_2\oplus \mathbb{Z}_2}$ and $\overline{X}^*_{\mathbb{Z}_2\oplus \mathbb{Z}_2}$  can be found in Figure \ref{figure:exampleD4}. The Hasse diagram of $X_{\mathbb{Z}_2\oplus \mathbb{Z}_2}$ is in black. The Hasse diagram of $\overline{X}_{\mathbb{Z}_2\oplus \mathbb{Z}_2}$ is in black, blue and red. The Hasse diagram of $\overline{X}^*_{\mathbb{Z}_2\oplus \mathbb{Z}_2}$ corresponds to the entire diagram, where we have in purple the new part added regard to $\overline{X}_{\mathbb{Z}_2\oplus \mathbb{Z}_2}$. 

\begin{figure}[h]
  \centering
 \includegraphics[scale=0.8]{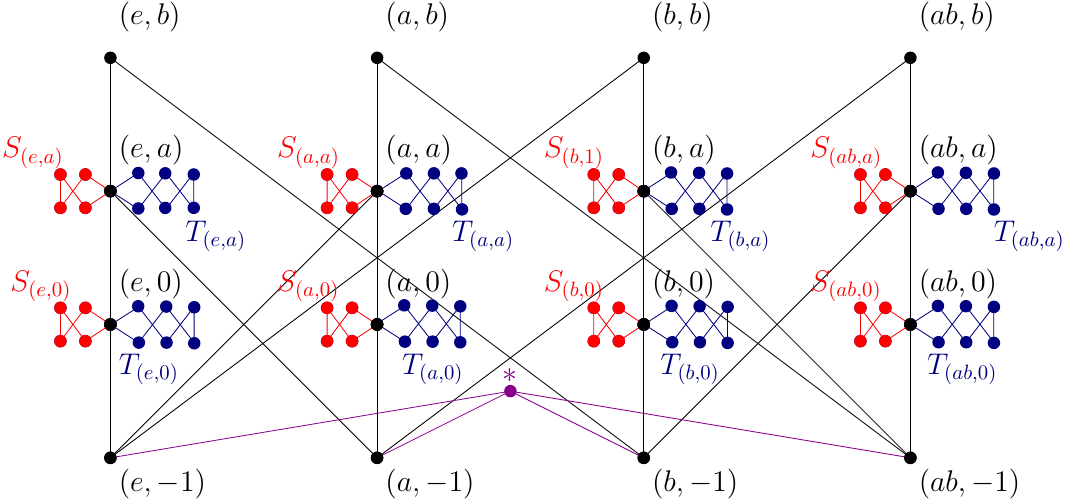}
     \caption{Hasse diagrams of $X_{\mathbb{Z}_2\oplus \mathbb{Z}_2}$,$\overline{X}_{\mathbb{Z}_2\oplus \mathbb{Z}_2}$ and $\overline{X}^*_{\mathbb{Z}_2\oplus \mathbb{Z}_2}$.}
     \label{figure:exampleD4}         
\end{figure}

\end{example}

\begin{lem}\label{lem:SyTsonasimetricos} Given  $g\in G$ and $0\leq \beta <max(S)$, $S_{(g,\beta)}\cup T_{(g,\beta)}$ as a subspace of $\overline{X}_G$ satisfies that its group of homeomorphisms is trivial.
\begin{proof}
We can deduce that $(g,\beta)$ is a fixed point for every homeomorphism since $(g,\beta)$ is not a maximal or minimal point. From here, it is easy to obtain the result.
\end{proof}
\end{lem}

\section{The group of homeomorphisms of $X_G$}\label{section:Automorphism(X_G)}
In this section, we prove Theorem \ref{thm:Aut(X_G)isoG}, the proof is essentially the same given by J.A Barmak and E.G. Minian   in \cite{barmak2009automorphism}, where they proved that every finite group $G$ can be realized as the group of homeomorphisms of a finite $T_0$ topological space (poset) with $n(r+2)$ points, where $|G|=n$ and $|S'|=r$. This result improves the results obtained by G. Birkhoff \cite{birkhoff1936order} and M.C. Thornton \cite{thornton1972spaces}, that used $n(n+1)$ and $n(2r+1)$ points respectively.

Before giving the proof of Theorem \ref{thm:Aut(X_G)isoG}, we study a property of $X_G$ that will be used in subsequent sections.

\begin{prop}\label{prop:preserveLevels} If $f:X_G\rightarrow X_G$ is a homeomorphism, then $f(G\times \{\beta\})=G\times \{\beta\}$ for every $\beta \in S$.
\begin{proof}

We consider $A=X_G\setminus \{(g,-1)\in X_G |g\in G\}$. $A$ is clearly the product of a discrete poset with a well-ordered set. Firstly, we show that $f(A)=A$. Then, we need to prove that $f(G\times \{ -1\})=G\times \{ -1\}$. We argue by contradiction,  suppose that there exists $(g,-1)\in X_G$ with $f(g,-1)=(h,\beta)$ for some $h\in G$ and $\beta>-1$. We have $f^{-1}(h,\beta)=(g,-1)$, so $f^{-1}(h,-1)\leq f^{-1}(h,\beta)$ since $(h,-1)<(h,\beta)$. By the minimality of $(g,-1)$, $f^{-1}(h,-1)=(g,-1)$, which leads to a contradiction with the injectivity of $f^{-1}$. Then, we can deduce that $f$ and $f^{-1}$ satisfy $f(G\times\{-1 \})=G\times\{-1 \}$ and $f^{-1}(G\times\{-1 \})=G\times\{-1 \}$. Thus, the restriction $f_{|A}:A\rightarrow A$ is a homeomorphism. $A$ has $|G|$ connected components, which are $A_g=\{(g,\beta)|\beta \geq 0\}$ for $g\in G$. If $(g,\beta)\in A$ satisfies that $f(g,\beta)=(h,\gamma)$ for some $h\in G$ and $\gamma\geq 0$, we get by the continuity of $f$ that $f(A_g)\subseteq A_h$. If there exists $y\in A_h\setminus  f(A_g)$, we can find an element $x\in f(A_g)$ such that $x<y$ or $x>y$, this can be done because $A_h$ is a well-ordered set. Hence, we know that $f^{-1}(x)\in A_g$. By the continuity of $f^{-1}$, we get that $y\in f(A_g)$, which leads to a contradiction. Therefore, $f(A_g)=A_h$. The only homeomorphism from a well-ordered set to itself is the identity, so $\gamma=\beta$. From here, we can deduce the desired result. 
\end{proof}
\end{prop}

\begin{proof}[Proof of Theorem \ref{thm:Aut(X_G)isoG}]
Given a group $G$, we consider the Alexandroff space $X_G$ constructed in Section \ref{section:constructionX_G}. We define $\varphi:G\rightarrow Aut(X_G)$ by $\varphi(g)(s,\beta)=(gs,\beta)$, where $(s,\beta)\in X_G$. We need to show that $\varphi$ is an isomorphism of groups. First of all, we check that $\varphi$ is well-defined. $\varphi(g):X_G\rightarrow X_G$ is clearly continuous because it preserves the order of $X_G$. By construction, $\varphi(g)$ is also bijective. The inverse of $\varphi(g)$ is $\varphi(g^{-1})$, which is also continuous. It is straightforward to check that $\varphi$ is a homomorphism of groups. 

We prove that $\varphi$ is a monomorphism of groups. Suppose that $\varphi(g)=id$, where $id:X_G\rightarrow X_G$ denotes the identity. Then, $(ge,-1)=\varphi(g)(e,-1)=(e,-1)$, where $e$ denotes the identity element of the group $G$, which implies that $g=e$.

Now, we verify that $\varphi$ is an epimorphism of groups. We consider  $f\in Aut(X_G)$. By Proposition \ref{prop:preserveLevels}, $f(e,-1)=(h,-1)$ for some $h\in G$. We also have that $\varphi(h)(e,-1)=(h,-1)$. We consider $Y=\{x\in X_G | f(x)=\varphi(h)(x) \}$. If we prove that $Y$ is open, $Y$ is closed and $X_G$ is a connected space, we get that $Y=X_G$ because $Y$ is non-empty ($(e,-1)\in Y$).

We prove that $Y$ is open. If $(g,\beta)\in Y$, we have $f_{|U_{(g,\beta)}},\varphi(h)_{|U_{(g,\beta)}}:U_{(g,\beta)}\rightarrow U_{f(g,\beta)}$ and $f(g,\beta)=(hg,\beta)=\varphi(h)(g,\beta)$. On the other hand, there is only one element for each $\gamma$ with $0\leq \gamma\leq \beta$ in $U_{(g,\beta)}$ and $U_{f(g,\beta)}$, $(g,\gamma)$ and $(hg,\gamma)$ respectively. Concretely, $U_{(g,\beta)}$ consists of points $(g,\gamma)$ with $-1\leq \gamma \leq \beta$ and points of the form $(g\gamma,-1)$ with $0<\gamma\leq \beta$, the description of $U_{f(g,\beta)}$ is similar. Hence, by Proposition \ref{prop:preserveLevels}, we can deduce that $f(g,\gamma)=(hg,\gamma)=\varphi(h)(g,\gamma)$ for every $0\leq \gamma\leq \beta$. It only remains to show that $f(g,-1)=\varphi(h)(g,-1)$ and $f(g\gamma,-1)=\varphi(h)(g\gamma,-1)$ for every $0<\gamma\leq \beta$. By the construction of $X_G$, $(g\gamma,-1)\prec (g,\gamma)$ for every $0< \gamma$. By Proposition \ref{prop:sucesores}, $f(g\gamma,-1)\prec f(g,\gamma)$. On the other hand, if $x\in U_{f(g,\beta)}$ satisfies $x\prec f(g,\gamma)=(hg,\gamma)$, we get that $x$ can only be of the form $(hg,\alpha)$ for some $0\leq \alpha<\gamma$ or $(hg\gamma,-1)$. Suppose that $f(g\gamma,-1)=(hg,\alpha)$ for some $0\leq \alpha<\gamma$, we get a contradiction with Proposition \ref{prop:preserveLevels}. Therefore, $f(g\gamma,-1)=(hg\gamma,-1)=\varphi(h)(g\gamma,-1)$. Finally, $(g,-1)\prec (g,0)$. By Proposition \ref{prop:sucesores}, $f(g,-1)\prec f(g,0)=(hg,0)$. Then, $f(g,-1)=(hg,-1)=\varphi(h)(g,-1)$. Thus, $\varphi(h)(x)=f(x)$ for every $x\in U_{(g,\beta)}$, that is to say, $U_{(g,\beta )}\subset Y$.


We prove that $Y$ is closed. We consider $(k,\beta)\in X_G\setminus Y$. By Proposition \ref{prop:preserveLevels}, $f(k,\beta)=(g,\beta)=\varphi(gk^{-1})(k,\beta)$ for some $g\in G$. Furthermore, $g\neq hk$ because otherwise we would get 
$$f(k,\beta)=(g,\beta)=(hk,\beta)=\varphi(hkk^{-1})(k,\beta)=\varphi(h)(k,\beta), $$
which leads to a contradiction with $(k,\beta)\in X_G\setminus Y$. We can repeat the same argument used before to get that $f_{|U_{(k,\beta)}}=\varphi(gk^{-1})_{|U_{(k,\beta)}}$, but $gk^{-1} \neq h$, so $f(y)=\varphi(gk^{-1}) (y)\neq \varphi(h)(y)$ for every $y\in U_{(k,\beta)}$, i.e.,  $U_{(k,\beta)}\cap Y=\emptyset$.

We prove that $X_G$ is connected showing that it is path-connected. We need to show that for every $x,y\in X_G$ there is a sequence $x=x_0,x_1,...,x_n=y$ with $x_i$ comparable to $x_{i+1}$ for every $i=0,...,n-1$. It is only necessary to check the existence of that sequence between points of the form $(g,-1),(h,-1)$ with $g,h\in G$ and $g\neq h$, the reason is the first relation of the partial order given in $X_G$. We prove that for every $g\in G$ there is a sequence of comparable points from $(g,-1)$ to $(e,-1)$. $S'$ is a generating set, so $g=\alpha_{1}^{d_1}\alpha_{2}^{d_2}...\alpha_{m}^{d_m}$, where $\alpha_{j}\in S'$ and $d_j=1$ or $-1$ with $j=1,...,m$. If $d_m=1$, $$(\alpha_{1}^{d_1}\alpha_{2}^{d_2}...\alpha_{m},-1)\prec (\alpha_{1}^{d_1}\alpha_{2}^{d_2}...\alpha_{{m-1}}^{d_{m-1}},\alpha_{m})>(\alpha_{1}^{d_1}\alpha_{2}^{d_2}...\alpha_{{m-1}}^{d_{m-1}},-1).$$

If $d_m=-1$, 
$$(\alpha_{1}^{d_1}\alpha_{2}^{d_2}...\alpha_{m}^{-1},-1)<(\alpha_{1}^{d_1}\alpha_{2}^{d_2}...\alpha_{m}^{-1},\alpha_{m})\succ (\alpha_{1}^{d_1}\alpha_{2}^{d_2}...\alpha_{{m-1}}^{d_{m-1}},-1). $$  

We need to combine these two steps inductively so as to obtain a sequence of comparable points from $(g,-1)$ to $(e,-1)$. From here, it can be deduced that for every $g,h\in G$ with $g\neq h$ there is a sequence of comparable points from $(g,-1)$ to $(h,-1)$. 

\end{proof}

\begin{remark}\label{remark:f=gsif(x)=g(x)} By the proof of Theorem \ref{thm:Aut(X_G)isoG}, the following property can be deduced: if $f,g\in Aut(X_G)$ and there exists $x\in X_G$ satisfying $f(x)=g(x)$, then $f=g$. 
\end{remark}

\section{The group of self-homotopy equivalences of $\overline{X}_G$ and $\overline{X}_G^*$}\label{section:selfhomotopyequivalences}
In this section, we show that every group $G$ can be realized as the group of (pointed) homotopy classes of (pointed) self-homotopy equivalences  of infinitely many non-homotopy-equivalent Alexandroff spaces, i.e.,  Theorem \ref{thm:equivInfinite} and Theorem \ref{cor:pointedHomotopy}. Given a non-trivial group $G$, we prove that $\mathcal{E}(\overline{X}_G)$ is isomorphic to $G$, we divide the proof into three lemmas.
\begin{lem}\label{lem:Aut(X_G)Aut(Xbarra_G)} Given a non-trivial group $G$,  $Aut(X_G)$ is isomorphic to $Aut(\overline{X}_G)$.
\begin{proof}
Firstly, we prove that every $f\in Aut(\overline{X}_G)$ satisfies $f(X_G)=X_G$, that is to say, $f_{|X_G}\in Aut (X_G)$. If $f\in Aut(\overline{X}_G)$, we have trivially that $f(M)=M$ and $f(N)=N$, where $M$ and $N$ denote the sets of maximal and minimal elements of $\overline{X}_G$. Then, we get $f(\overline{X}_G\setminus (M\cup N))=\overline{X}_G\setminus (M\cup N)$. 

Now, we prove that $f(G\times\{-1\})=G\times\{-1\}$. We argue by contradiction, suppose there exists $g\in G$ such that $f(g,-1)\in N\setminus(G\times \{ -1\})$. If $f(g,-1)=x$ with $x\in \{C_{(h,\alpha)},D_{(h,\alpha)},H_{(h,\alpha)},I_{(h,\alpha)},J_{(h,\alpha)} \}$ for some $h\in G$ and $0\leq \alpha< max(S)$, the set $Z=\{y\in \overline{X}_G|x\prec y \}$ has cardinality two and it contains at least one element in $M$. The set $Y=\{y\in \overline{X}_G|(g,-1)\prec y \}$ has cardinality at least two. If $|S'|\geq 2$,  $(g,0)$ and $(g\beta^{-1},\beta)$ with $\beta < max(S)$ are points in $Y\setminus (M\cup N)$. By Proposition \ref{prop:sucesores} and the injectivity of $f$, $f(Y\setminus (M\cup N))=Z$, but it contradicts $f(\overline{X}_G\setminus (M\cup N))=\overline{X}_G\setminus (M\cup N)$. If $|S'|=1$, we know that $f(G\times \{0 \})=G\times \{0\}$ since the points of the form $(k,0)$ for some $k\in G$ are the only points that are not maximal or minimal points. If $f(g,-1)=x$ with $x\in \{C_{(h,0)},I_{(h,0)},J_{(h,0)} \}$ for some $h\in G$, we get a contradiction studying the image by $f$ of $(g,-1)\prec (g,0)\prec (g,max(S))$. If $f(g,-1)=x$ with $x\in \{D_{(h,0)},H_{(h,0)} \}$ for some $h\in G$,  we get a contradiction studying the image by $f$ of $(g,-1)\prec (g\ max(S)^{-1},max(S))\succ (g\ max(S)^{-1},0)$. We have shown that $f(G\times\{-1\})\subseteq G\times\{-1\}$. The equality follows due to the fact that $f$ is a homeomorphism.


If $S$ is finite or $S$ is not finite and has a maximum such that $max(S)$ has an immediate predecessor, we prove that $f(G\times \{max(S) \})=G\times \{max(S)\}$. We argue by contradiction, suppose there exists $g\in G$ such that $f(g,max(S))=x$ with $x\in \{A_{(h,\alpha)},$ $B_{(h,\alpha)},$ $E_{(h,\alpha)},$ $F_{(h,\alpha)},$ $G_{(h,\alpha)}\}$ for some $h\in G$ and $0\leq \alpha< max(S)$. We repeat the previous argument. The set $Z=\{y\in \overline{X}_G|x\succ y \}$ has cardinality two and it is a subset of $\overline{X}_G\setminus (G\times \{ -1\})$. The set $Y=\{y\in \overline{X}_G|(g,max(S))\succ y \}$ has cardinality two, its elements are $\{(g,\alpha)$ and $(g$ $ max(S),-1)$ for some $0\leq\alpha<max(S)$. By Proposition \ref{prop:sucesores} and the injectivity of $f$, $f(Y)= Z$, but it contradicts $f(G\times\{-1\})=G\times\{-1\}$. Again, $f(G\times \{max(S) \})=G\times \{max(S)\}$ follows from the fact that $f$ is a homeomorphism.

By the previous arguments, it can be deduced that $f(X_G)=X_G$. Moreover, by Proposition \ref{prop:preserveLevels}, $f(G\times \{\beta\})=G\times \{ \beta\}$ for every $\beta \in S$. By Lemma \ref{lem:SyTsonasimetricos} and Proposition \ref{prop:sucesores}, for every $(g,\beta)$ with $g\in G$ and $0\leq \beta < max(S) $, we get  $f(S_{(g,\beta)})=S_{f(g,\beta)}$ and $f(T_{(g,\beta)})=T_{f(g,\beta)}$.

We define $\varphi:Aut(\overline{X}_G)\rightarrow Aut(X_G)$ given by $\varphi(f)=f_{|X_G}$. It is immediate to show that $\varphi$ is a homomorphism of groups and is well-defined. We prove that $\varphi$ is a monomorphism of groups. If $f,h\in Aut(\overline{X}_G)$ with $f\neq h$, it means that there exists a point $x\in \overline{X}_G$ such that $f(x)\neq h(x)$. Suppose $x\in (S_{(g,\beta)}\cup T_{(g,\beta)})\setminus \{ (g,\beta)\}$ for some $(g,\beta)\in X_G$, by the previous paragraph, it can be deduced that $f(g,\beta)\neq h(g,\beta)$. Hence, it can be found $y\in X_G$ such that $f_{|X_G}(y)\neq h_{|X_G}(y)$. Now, we prove that $\varphi$ is an epimorphism of groups. If $f\in Aut(X_G)$, we extend $f:X_G\rightarrow X_G$ to $f':\overline{X}_{G}\rightarrow \overline{X}_G$ just declaring that $f'(x)=f(x)$ for every $x\in X_G$ and $f'(S_{(g,\beta)})=S_{f(g,\beta)}$, $f'(T_{(g,\beta)})=T_{f(g,\beta)}$ for every $g\in G$ and $0\leq\beta<max(S)$. It is easy to check that $f'\in Aut(\overline{X}_G)$. Then, $\varphi(f')=f$.
\end{proof}
\end{lem}

\begin{lem}\label{lem:simplificar} If $X$ is an Alexandroff space, $a\in X$ is a maximal (resp. minimal) element such that $U_a\setminus \{a\}$ (resp. $F_a\setminus\{ a\}$) is not connected and $r:X\rightarrow X$ is a comparative retraction, then $r(a)=a$.
\begin{proof}
We argue by contradiction, suppose that $r(a)<a$, the case $r(a)>a$ is not possible since $a$ is a maximal element. By hypothesis, $U_a\setminus\{a \}=V \cup W$, where $V$ and $W$ are disjoint non-empty open sets. We can suppose that $x=r(a)\in V$, then $U_x\subseteq V$. We take $y\in W$, so $U_y\subseteq W$. By the continuity of $r$, we get $r(y)\leq r(a)$. We also have that $r$ is a comparative retraction, so $r(y)\leq y$ or $r(y)\geq y$. Suppose that $r(y)\leq y$, we also have that $r(y)\leq x$, therefore, $r(y)\in U_x\cap U_y\subseteq V\cap W $, which entails a contradiction. If $r(y)\geq y$, we have $x\geq r(y)\geq y$ and $y\in U_x$, which leads to a contradiction.


The case when $a$ is a minimal element follows from the previous case. If $r:X\rightarrow X$ is a comparative retraction, it is also a comparative retraction when we consider the opposite order in $X$. With the opposite order we are in the previous conditions, that is to say, $a$ is now a maximal element and $U_{a}\setminus \{a \}$ is not connected. Therefore, we get that $r(a)=a$.
\end{proof}
\end{lem}

\begin{lem}\label{lem:E(Xbarra_G)Aut(Xbarra_G)}
Given a non-trivial group $G$, $\mathcal{E}(\overline{X}_G)$ is isomorphic to $Aut(\overline{X}_G)$. 
\begin{proof}
We need to prove that $\overline{X}_G$ is locally a core. Then, applying Lemma \ref{lem:clave}, we get the desired result. Firstly, we show that $\overline{X}_G$ is a $\mathcal{C}$-core. Hence, we need to verify that the only comparative retraction of $X$ is the identity. 

We take a comparative retraction $r$. If $x\in (S_{(g,\beta)} \cup T_{(g,\beta)})\setminus \{(g,\beta)\}$ or $x\in $ $\{(g,-1),$ $(g,max(S))\}$, where $g\in G$ and $0\leq \beta <max(S)$,  we are trivially in the hypothesis of Lemma \ref{lem:simplificar}, so $r(x)=x$.

It only remains to study the points of the form $(g,\beta)$ with $0\leq \beta < max(S)$ and $g\in G$. We argue by contradiction. Suppose that $r(g,\beta)\neq (g,\beta)$. Then, $r(g,\beta)<(g,\beta)$ or $r(g,\beta)>(g,\beta)$. Suppose $r(g,\beta)<(g,\beta)$. If $r(g,\beta)\neq D_{(g,\beta)}$, we get $D_{(g,\beta)}=r(D_{(g,\beta)})<r(g,\beta)<(g,\beta)$, which leads to a contradiction because $D_{(g,\beta)}\prec (g,\beta)$. If $r(g,\beta)=D_{(g,\beta)}$, we get a contradiction because we get $H_{(g,\beta)}=r(H_{(g,\beta)})\leq r(g,\beta)=D_{(g,\beta)}$. Suppose $r(g,\beta)>(g,\beta)$. We can repeat the same arguments used before. If $r(g,\beta)\neq B_{(g,\beta)}$, we have $B_{(g,\beta)}=r(B_{(g,\beta)})>r(g,\beta)>(g,\beta)$, which leads to a contradiction with $B_{(g,\beta)}\succ (g,\beta)$. If $r(g,\beta)= B_{(g,\beta)}$, we get $E_{(g,\beta)}=r(E_{(g,\beta)})\geq r(g,\beta)=B_{(g,\beta)}$, which entails a contradiction.

We have shown that $\overline{X}_G$ is a $\mathcal{C}$-core. Now, we  prove that $\overline{X}_G$ is locally a core. For every $x\in S_{(g,\beta)}\cup  T_{(g,\beta)}$, where $0\leq \beta <max(S)$ and $g\in G$, we consider $A_x=S_{(g,\beta)}\cup T_{(g,\beta)}$. Suppose $x=(g,-1)$ with $g\in G$. Then, $(g,-1)\prec (g\gamma^{-1},\gamma)$ for some $\gamma\in S'$. If $\gamma\neq max(S)$, we consider $A_x=S_{(g,0)}\cup T_{(g,0)}\cup S_{(g\gamma^{-1},\gamma)}\cup T_{(g\gamma^{-1},\gamma)}\cup \{(g,-1)\} $. If $S$ is finite or $S$ is not finite but has a maximum such that $max(S)$ has an immediate predecessor and $\gamma=max(S)$, we consider $A_x=S_{(g,0)}\cup T_{(g,0)}\cup S_{(g\gamma^{-1},\alpha)}\cup T_{(g\gamma^{-1},\alpha)}\cup \{(g,-1),(g\gamma^{-1},\gamma)\}$, where $(g\gamma^{-1},\alpha)\prec (g\gamma^{-1},\gamma)$. If $S$ is finite or $S$ is not finite but has a maximum such that $max(S')$ has an immediate predecessor and $x=(g,max(S))$, there exists $\alpha\in S$ such that $(g,\alpha)\prec (g,max(S))$ and we consider $A_x= S_{(g,\alpha)}\cup  T_{(g,\alpha)}\cup S_{(g \ max(S),0)}\cup T_{(g\ max(S),0)} \cup \{x,(g$ $max(S),-1)\} $. It is immediate to show that $A_x$ satisfies the property asked in Definition \ref{def:locallycore} for every $x\in \overline{X}_G$.

\end{proof}
\end{lem}

\begin{remark}\label{rem:pointedCasenoVabien} By the proofs of Remark \ref{remark:f=gsif(x)=g(x)}, Lemma \ref{lem:Aut(X_G)Aut(Xbarra_G)} and Lemma \ref{lem:E(Xbarra_G)Aut(Xbarra_G)}, it can be deduced that the set $L_{x,y}=\{ f:(\overline{X}_G,x)\rightarrow (\overline{X}_G,y)| f(x)=y $ and $f\in \mathcal{E}(\overline{X}_G) \}$ has cardinality at most $1$. 
\end{remark}

In general, for an arbitrary Alexandroff space we cannot expect to obtain an isomorphism of groups between the group of homeomorphisms and the group of homotopy classes of self-homotopy equivalences.
\begin{example}\label{example:autnoequiv}
We consider $A=\{a,b,c,d,e \}$, where we use the topology associated to the following partial order: $a,b<d,c,e$ and $c<e$. We study $Aut(A)$. A homeomorphism $f:A\rightarrow A$ preserves the order and therefore should send maximal chains to maximal chains. In $A$, there are two maximal chains of three elements, $a<c<e$ and $b<c<e$. From here, it is easy to deduce that $e,d$ and $c$ are fixed points for every homeomorphism. Then, $Aut(A)\simeq \mathbb{Z}_2$. On the other hand, $A^c=\{a,b,c,d\}$ is the core of $A$ because $e$ is clearly a down beat point and $A^c$ does not contain beat points. Hence, $\mathcal{E}(A)\simeq \mathcal{E}(A^c)$. In addition, by Lemma \ref{lem:clave}, we get $Aut(A^c)\simeq \mathcal{E}(A^c)$. From here, it is immediate that $\mathcal{E}(A^c)$ is the Klein four-group.  We describe the two generators $f$ and $g$ of $\mathbb{Z}_2\oplus \mathbb{Z}_2\simeq \mathcal{E}(A)$. $f$ is given by $f(a)=b, f(b)=a, f(c)=c, f(d)=d$ and $g$ is given by $g(c)=d,g(d)=c,g(a)=a,g(b)=b$. A schematic situation in the Hasse diagrams can be seen in Figure \ref{figure:ejemploAutnoEquiv}.
\begin{figure}[h]
  \centering
    \includegraphics[scale=0.6]{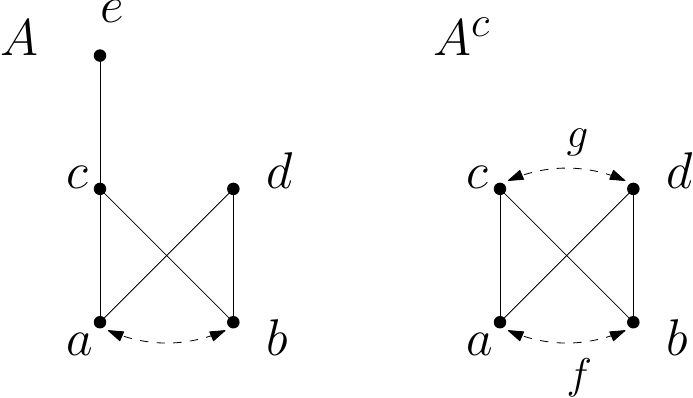}
     \caption{Hasse diagrams of $A$ and $A^c$.}
     \label{figure:ejemploAutnoEquiv}
\end{figure}
\end{example}

\begin{theorem}\label{thm:GisoE(X)1}
Every group can be realized as the group of self-homotopy equivalences of an Alexandroff space.
\begin{proof}
Given a non-trivial group $G$, we consider $X_G$ and $\overline{X}_G$. By the proof of Theorem \ref{thm:Aut(X_G)isoG}, $G\simeq Aut(X_G)$. By Lemma \ref{lem:Aut(X_G)Aut(Xbarra_G)} and Lemma \ref{lem:E(Xbarra_G)Aut(Xbarra_G)}, we get $\mathcal{E}(\overline{X}_G)\simeq Aut(X_G)$. If $G$ is the trivial group, we only need to consider the Alexandroff space given by one point.
\end{proof}
\end{theorem}
%

\begin{figure}[h]
  \centering
    \includegraphics[scale=0.7]{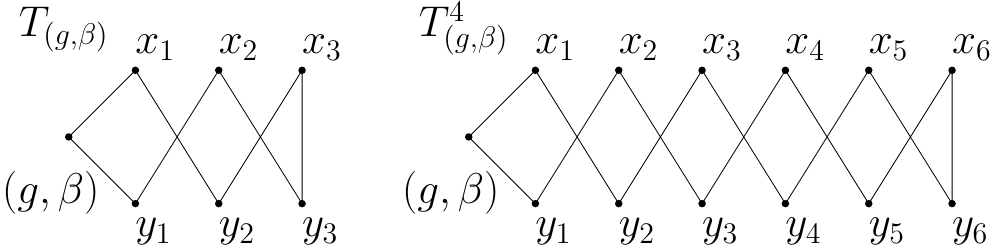}
     \caption{Hasse diagrams of $T_{(g,\beta)}^1$ and $T^4_{(g,\beta)}$.}
     \label{figure:TnyT}
\end{figure}

Given a group $G$, a slight modification of the construction made in Section \ref{section:constructionX_G} can provide infinitely many Alexandroff spaces satisfying that their groups of homotopy classes of self-homotopy equivalences are isomorphic to $G$. We only need to change $T_{(g,\beta)}$ by $T_{(g,\beta)}^n$ for every $(g,\beta)\in \overline{X}_G$, where $0\leq \beta <max(S)$ and $n \in \mathbb{N}$, so as to get $\overline{X}_G^n$. $T_{(g,\beta)}^n$ consists of $2n+5$ points, concretely, $T_{(g,\beta)}^n=\{x_1,...,x_{2+n},y_1,...,y_{2+n}, (g,\beta) \}$, where $x_i$ denotes a maximal element and $y_i$ denotes a minimal element for $i=1,...,2+n$. The relations are given by the following formulas:
\begin{align}\label{eq:relaciones}
(g,\beta)<x_1>y_2<x_3>...<x_{1+n}>y_{2+n}<x_{2+n}>y_{1+n}<x_{n}>...<x_2>y_1<(g,\beta), 
\end{align}
\begin{align}\label{eq:relaciones2}
(g,\beta)<x_1>y_2<x_3>...>y_{1+n}<x_{2+n}>y_{2+n}<x_{1+n}<y_{n}>...<x_2>y_1<(g,\beta),
\end{align}
where the first case is considered when $n$ is even and the second case when $n$ is odd. An example of the Hasse diagrams can be seen in Figure \ref{figure:TnyT}.
It is clear that $T_{(g,\beta)}^1=T_{(g,\beta)}$ and $\overline{X}_G^1=\overline{X}_G$. We consider 
$$\overline{X}_G^n=X_G\cup (\bigcup_{\substack{{(g,\beta)\in G\times S} \\ 0\leq \beta <max(S)}}(S_{(g,\beta)}\cup T^n_{(g,\beta)})),$$
where we extend the partial order of $X_G\cup(\bigcup_{\substack{{(g,\beta)\in G\times S} \\ 0\leq \beta <max(S)}}S_{(g,\beta)})$ as  a subspace of $\overline{X}_G$ to $\overline{X}^n_G$. If $x_1,y_1\in T_{(g,\beta)}^n$, we identify them with $E_{(g,\beta)}$ and $H_{(g,\beta)}$ respectively. Then, we consider the relations 3. 5. 8. and 9. of $\overline{X}_G$. Finally, we consider the relations given in Formula (\ref{eq:relaciones}) and (\ref{eq:relaciones2}).
We need to combine the previous results and constructions to obtain the proof of  Theorem \ref{thm:equivInfinite}.
\begin{proof}[Proof of Theorem \ref{thm:equivInfinite}]
Given a non-trivial group $G$, we consider $\{\overline{X}_G^n\}_{n\in \mathbb{N}}$. We can prove that $Aut(\overline{X}_G^n)\simeq Aut(X_G)$ using similar arguments than the ones used in the proof of Lemma \ref{lem:Aut(X_G)Aut(Xbarra_G)}. Following the same techniques of the proof of Lemma \ref{lem:E(Xbarra_G)Aut(Xbarra_G)}, it can be obtained that $\mathcal{E}(\overline{X}_G^n)\simeq Aut(\overline{X}_G^n)$. By the proof of Theorem \ref{thm:Aut(X_G)isoG}, we get $Aut(X_G)\simeq G$. 

If $\overline{X}^n_G$ is homotopy equivalent to $\overline{X}^m_G$, where $m\neq n$, there exist two continuous functions $f:\overline{X}^n_G\rightarrow \overline{X}^m_G $, $g:\overline{X}^m_G \rightarrow \overline{X}^n_G$ such that $f\circ g\simeq id_{\overline{X}^m_G}$ and $g\circ f\simeq id_{\overline{X}^n_G}$. Without loss of generality we assume that $n>m$. By Theorem \ref{thm:kukielaequivIdentity}, $f\circ g= id_{\overline{X}^m_G}$ and  $g\circ f= id_{\overline{X}^n_G}$. Then, $f$ is a homeomorphism and $g$ is the inverse of $f$. Therefore, $f(G\times \{ \beta \})=G\times \{ \beta \}$ with $0\leq \beta <max(S)$. To prove the last assertion we need to argue as we did in the proof of Lemma \ref{lem:Aut(X_G)Aut(Xbarra_G)}. If $M^s$ and $N^s$ denote the set of maximal and minimal points of $\overline{X}^s_G$, where $s=n,m$, we have $f(\overline{X}^n_G\setminus (M^n\cup N^n))=\overline{X}^m_G\setminus (M^m\cup N^m)$, where $\overline{X}^s_G \setminus (M^s\cup N^s)$ is the product of a discrete poset with a well-ordered set for $s=n,m$. By the proof of Proposition \ref{prop:preserveLevels}, it can be deduced the desired assertion. By Proposition \ref{prop:sucesores} and the previous assertion, for every $(g, \beta)$, where $g\in G$ and $0\leq \beta<max(S)$, we can deduce that $f(S_{(g,\beta )}\cup T^n_{(g,\beta)})\subseteq S_{(g,\beta )}\cup T^m_{(g,\beta)}$, which leads to a contradiction with the injectivity of $f$.

If $G$ is the trivial group, we only need to consider $\{S_{(e,0)}\cup T^n_{(e,0)}\}_{n\in \mathbb{N}}$. We can apply the same techniques and a generalization of Lemma \ref{lem:SyTsonasimetricos} so as to prove the desired result.

\end{proof}
\begin{remark} We can deduce from the proof of Theorem \ref{thm:equivInfinite} an analogue of that theorem for the group of homeomorphisms, i.e., every group $G$ can be realized as the group of homeomorphisms of infinitely many (non-homeomorphic) Alexandroff spaces. 
\end{remark}

\begin{proof}[Proof of Theorem \ref{cor:pointedHomotopy}]
Given a non-trivial group $G$, the Alexandroff spaces $\{\overline{X}_G^n\}_{n\in \mathbb{N}}$ used in the proof of Theorem \ref{thm:equivInfinite} are far from satisfying that their group of pointed homotopy classes of pointed self-homotopy equivalences is isomorphic to $G$. The reason for the previous fact is Remark \ref{rem:pointedCasenoVabien}, it can be adapted to $\overline{X}_G^n$ for every $n\in \mathbb{N}$. Nevertheless, it will be only necessary to add one extra point to $\overline{X}_G^n$ so as to obtain the desired result. We consider $\overline{X}_G^{n*}=\overline{X}_G^{n}\cup \{*\}$, where we extend the partial order of $\overline{X}_G^{n}$ to $\overline{X}_G^{n*}$ declaring that $*>(g,-1)$ for every $g\in G$. Firstly, we prove that $\overline{X}_G^{n*}$ is a $\mathcal{C}$-core. We need to verify that the only comparative retraction $r$ is the identity. If  $x\in (S_{(g,\beta)}\cup T^n_{(g,\beta)} )\setminus \{(g,\beta)\}$ or $x\in \{ (g,-1),*,(g,max(S)) \}$, where $0\leq \beta<max(S)$ and $g\in G$, then $x$ satisfies the hypothesis of Lemma \ref{lem:simplificar}. Hence, we need to show that $r(x)=x$ for the points of the form $(g,\beta)$, where $0\leq \beta <max(S)$ and $g\in G$. We only need to repeat the same arguments we used in the proof of Theorem \ref{thm:equivInfinite} for that sort of points so as to conclude.


We prove that $\overline{X}_G^{n*}$ is locally a core. If $x\in S_{(g,\beta)}\cup  T_{(g,\beta )}^n $, where $g\in G$ and $0\leq \beta< max(S)$, we consider $A_x= S_{(g,\beta)} \cup T_{(g,\beta )}^n$. Suppose $x=(g,-1)$ with  $g\in G$. Then, $(g,-1)\prec (g\gamma^{-1},\gamma)$ for some $\gamma\in S' $. If $\gamma\neq max(S)$, we consider $A_x=S_{(g,0)}\cup T_{(g,0)}^n\cup S_{(g\gamma^{-1},\gamma)}\cup T_{(g\gamma^{-1},\gamma)}^n\cup \{(g,-1)\} $. If $S$ is finite or $S$ is not finite but has a maximum such that $max(S)$ has an immediate predecessor and $\gamma=max(S)$, we consider $A_x=S_{(g,0)}\cup T^n_{(g,0)}\cup S_{(g\gamma^{-1},\alpha)}\cup T^n_{(g\gamma^{-1},\alpha)}\cup \{ (g,-1),(g\gamma^{-1},\gamma)\}$, where $(g\gamma^{-1},\alpha)\prec (g\gamma^{-1},\gamma)$. If $S$ is finite or $S$ is not finite but has a maximum such that $max(S)$ has an immediate predecessor and $x=(g,max(S))$ for some $g\in G$, there exists $\alpha\in S$ such that $(g,\alpha)\prec (g,max(S))$ and we consider $A_x= S_{(g,\alpha)}\cup  T_{(g,\alpha)}^n\cup S_{(g \ max(S),0)}\cup T_{(g\ max(S),0)}^n \cup \{x,(g$ $max(S),-1)\} $. Finally, if $x=*$, we consider the sets $A_{(g,-1)}$ and $A_{(h,-1)}$ defined before, where $g,h\in G$ and $g\neq h$. Therefore, $A_*=\{*\} \cup A_{(g,-1)} \cup {A_{(h,-1)}}$. For every $x\in \overline{X}_G^n $, $A_x$ satisfies the property asked in Definition \ref{def:locallycore}. Thus, $\overline{X}_G^{n*}$ is locally a core. By Lemma \ref{lem:clave}, $\mathcal{E}(\overline{X}_G^{n*})$ is isomorphic to $Aut(\overline{X}_G^{n*})$.

We prove that $Aut(\overline{X}^n_G)$ is isomorphic to $ Aut(\overline{X}_G^{n*})$. To do that, we show that $*$ is a fixed point for every $f\in Aut(\overline{X}^{n*}_G)$. We argue by contradiction, suppose $*$ is not a fixed point for $f\in  Aut(\overline{X}^{n*}_G)$, we study cases. By construction, $* $ is a maximal element. Then, $f(*)=x$ with $x$ being a maximal element of $\overline{X}^n_G$, that is to say,  $x$
is a maximal element in $S_{(g,\beta)}\cup T_{(g,\beta)}^n$ or is equal to $(g,max(S))$, where $g\in G$ and $0\leq \beta <max(S)$. We consider $Y_{z}=\{ y\in \overline{X}^n_G| y\prec z\}$. Then, we have $|Y_{*}|=|G|$ and $|Y_{x}|=2$. By Proposition \ref{prop:sucesores}, $f(Y_*)\subseteq Y_x$. Then, if $|G|>2$, we get a contradiction with the injectivity of $f$. If $|G|=2$, $G$ is isomorphic to $\mathbb{Z}_2$. Since $*$ is a maximal element, $f(*)$ is a maximal element, we study cases. If $f(*)$ is a maximal element in $S_{(i,0)}\cup T_{(i,0)}^n$, where $i\in \{0,1 \}$, we get a contradiction studying the image by $f$ of $*\succ (0,-1)\prec(0,0)\prec(0,1)$. If $f(*)=(i,1)$, where $i\in \{0,1 \} $, we get a contradiction studying the image by $f^{-1}$ of $(i,1)\succ (i,0) \succ (i,-1)$. Thus, the only possibility is  $f(*)=*$, which means that $*$ is a fixed point for every $f\in Aut(\overline{X}_G^{n*})$. It follows that $Aut(\overline{X}^n_G)$ is isomorphic to $ Aut(\overline{X}_G^{n*})$.

Now, we consider the group of pointed homotopy classes of pointed self-homotopy equivalences of the pointed space $(\overline{X}_G^{n*},*)$. We also denote that group by $\mathcal{E}((\overline{X}_G^{n*},*))$. 
By the proof of Lemma \ref{lem:clave} and Remark \ref{rem:AutEqSonIguales}, it can be deduced that $*$ is a fixed point for every self-homotopy equivalence $f\in \mathcal{E}(\overline{X}_G^{n*})$. Then, we get $\mathcal{E}(\overline{X}_G^{n*})=\mathcal{E}((\overline{X}_G^{n*},*))$. By the proof of Theorem \ref{thm:equivInfinite}, we obtain $Aut(\overline{X}_G^n)\simeq G$. From here, we deduce the desired result.

Finally, we prove the following: $\overline{X}_G^{n*}$ is homotopy equivalent to $\overline{X}_G^{m*}$ if and only if  $m=n$.  We can argue as we did in the proof of Theorem \ref{thm:equivInfinite}. Suppose that $n>m$ and $\overline{X}_G^{n*}$ is homotopy equivalent to $\overline{X}_G^{m*}$. Therefore, it can be deduced that there exists a homeomorphism $f:\overline{X}_G^{n*}\rightarrow \overline{X}_G^{m*}$. We consider $M_*^s=M^s\cup \{ *\}$ and $N^s_*=N^s$, where $M^s$ and $N^s$ denote the sets of maximal and minimal elements of $\overline{X}_G^s$ with $s=n,m$. From here, we only need to repeat the same arguments used in the proof of Theorem \ref{thm:equivInfinite} so as to conclude.

If $G$ is the trivial group, we only need to consider $\{(S_{(e,0)}\cup T^n_{(e,0)},(e,0))\}_{n\in \mathbb{N}}$, where we have that $(e,0)$ is a fixed point for every homeomorphism. From here, it can be deduced the result.
\end{proof}

\begin{remark} If $G$ is a finite group and $S'$ is a non-trivial generating set of $G$, we can consider $\overline{X}'_G=X_G\cup (\bigcup_{g\in G, 0\leq \beta < max(S)} S_{(g,\beta)})$, where the partial order of $\overline{X}'_G$ is the one given as a subspace of $\overline{X}_G$. Since $G$ is finite, $\overline{X}'_G$ and $\overline{X}_G$ are finite topological spaces. Concretely,  $\overline{X}_G$ has $|G|(|S'|+2)+10|G||S'|$ points, while $\overline{X}'_G$ has $|G|(|S'|+2)+4|G||S'|$ points. We can repeat the same techniques used in the proof of Lemma \ref{lem:Aut(X_G)Aut(Xbarra_G)} to prove the analogue for $\overline{X}'_G$. Similarly, it can be obtained an analogue of Lemma \ref{lem:E(Xbarra_G)Aut(Xbarra_G)} for $\overline{X}_G'$, where the set $A_x$ of Definition \ref{def:locallycore} can be taken as $A_x=\overline{X}'_G$ for every $x\in \overline{X}_G'$. Hence, it can be proved that $Aut(\overline{X}_G')\simeq \mathcal{E}(\overline{X}_G')\simeq G$ for every finite group $G$. Analogously, the pointed case can also be obtained. Thus, for the finite case, it can be used an Alexandroff space with lower cardinality than in the general case.
\end{remark}


\section{Some properties of $\overline{X}_G$}\label{Section:McCord}
Let $G$ be a countable group, i.e., it is countable as a set or it has a countable set of generators. In this section, we study the weak homotopy type of the space $\overline{X}_G$. 

Let $S'$ be a set of non-trivial generators of $G$. If $S'$ is countable, we can denote the elements of $S'$ by $s_i$ with $i\in I\subseteq \mathbb{N}$. If $S'$ is finite, $I$ can be taken as $I=\{1,2,...,|S'| \}$. We declare $s_i<s_j$ if and only if $i<j$. $S'$ with the previous relation is a well-ordered set. In addition, for every $x\in S'$, $U_x$ is a finite set. Thus, taking $S'$ and the previous relation in the construction of $\overline{X}_G$, it can be deduced that $\overline{X}_G$ is locally finite.


We define the undirected graph $H_u(\overline{X}_G)$ given by the Hasse diagram of $\overline{X}_G$, the set of vertices are the points of $\overline{X}_G$ and there is an edge between two vertices $x$ and $y$ if and only if $x\prec y$ ($x\succ y$). $H_u(\overline{X}_G)$ can be seen as a one-dimensional CW-complex. 

\begin{prop}\label{prop:homotopyTypeX_G} If $G$ is a countable infinite group, then  $|\mathcal{K}(\overline{X}_G)|$ is homotopy equivalent to $\bigvee_{\mathbb{N}}S^1$.  If $G$ is a finite group, $|\mathcal{K}(\overline{X}_G)|$ is homotopy equivalent to the wedge sum of $3nr-n+1$ copies of $S^1$, where $|G|=n$ and $|S'|=r$. 
\begin{proof}
Firstly, we prove that $H_u(\overline{X}_G)$ and $ |\mathcal{K}(\overline{X}_G)|$ have the same homotopy type. The idea of the proof is to show that the natural inclusion $i:H_u(\overline{X}_G)\rightarrow |\mathcal{K}(\overline{X}_G) |$ is a weak homotopy equivalence between two CW-complexes. Then, by a well-known theorem of Whitehead, we would get that $H_u(\overline{X}_G)$ is homotopy equivalent to $| \mathcal{K}(\overline{X}_G)|$.

By the construction of the McCord complex, $H_u(\overline{X}_G)$ is a subcomplex of $\mathcal{K}(\overline{X}_G)$, so $i$ is well-defined and continuous. By Theorem \ref{thm:McCord}, we know that there is a weak homotopy equivalence $f:|\mathcal{K}(\overline{X}_G)|\rightarrow \overline{X}_G$. Furthermore, by the proof of Lemma 3 in \cite{mccord1966singular}, $f^{-1}(U_x)=\bigcup_{y\in U_x }star(y)$, where $star(y)$ denotes the union of all open simplices from $\mathcal{K}(\overline{X}_G)$ containing $y$ as a vertex. 
By the proof of Corollary 11 of Chapter 3 in \cite{spanier1981algebraic}, we can get that $|\mathcal{K}(U_x)|$ is a strong deformation retract of $f^{-1}(U_x)$, since it can be proved that $\mathcal{K}(\overline{X}_G)\setminus f^{-1}(U_x)$ is the largest subcomplex of $\mathcal{K}(\overline{X}_G)$ disjoint from $\mathcal{K}(U_x)$ and $\mathcal{K}(U_x)$ is a full subcomplex of $\mathcal{K}(\overline{X}_G)$. Finally, $U_x$ is contractible because it has a maximum. Then, it can be deduced that $|\mathcal{K}(U_x)|$ is also contractible.

 

$\mathcal{U}$ denotes the basis-like open cover given by $\{U_x\}_{x\in \overline{X}_G}$. It is straightforward to check that $f^{-1}(\mathcal{U})=\{f^{-1}(U_x)|U_x\in \mathcal{U} \}$ is a basis-like open cover for $|\mathcal{K}(\overline{X}_G)|$.

Now, we denote $B_x=\bigcup_{y\in U_x\subset H_u(U_x)} star(y)$, here,  $star(y)$ denotes the union of all open simplices from $H_u(\overline{X}_G)$ containing $y$ as a vertex. Then, we have that $i^{-1}(f^{-1}(U_x))=B_x$. It is trivial to show that $H_u(U_x)$ is a full subcomplex of $H_u(\overline{X}_G)$. We can repeat the previous arguments so as to prove that $H_u(U_x)$ is a strong deformation retract of $B_x$, i.e., we can use the proof of Corollary 11 of Chapter $3$ in \cite{spanier1981algebraic}. 

On the other hand, $H_u(U_x)$ is contractible for every $x\in \overline{X}_G$ since $H_u(U_x)$ is a tree. The vertices of $H_u(U_x)$ are the points of $\overline{X}_G$ that are smaller or equal to $x$. In Figure \ref{figure:grafosU}, we have the different graphs that can appear when we consider $H_u(U_x)$. If $x\in \{ A_{(g,\beta)},F_{(g,\beta)},G_{(g,\beta)}\}$, where $g\in G$ and $0\leq \beta <max(S)$, $H_u(U_x)$ is isomorphic to graph A). If $x\in\{ C_{(g,\beta)},$ $D_{(g,\beta)},$ $H_{(g,\beta)},$ $I_{(g,\beta)},$ $J_{(g,\beta)}$ $,(g,-1)\}$, where $g\in G$ and $0\leq \beta <max(S)$, $H_u(U_x)$ is isomorphic to graph B). If $x=(g,0)$, for some $g\in G$, $H_u(U_x)$ is isomorphic to graph C). If $x\in \{B_{(g,\beta)}, E_{(g,\beta)} \}$, where $g\in G$ and  $0\leq \beta <max(S)$, $H_u(U_x)$ is isomorphic to graph D), which is the graph in blue and black. If $x=(g,s_n)$ for some $g\in G$ and $n\in I\subseteq \mathbb{N}$, $H_u(U_x)$ is isomorphic to graph E). If $x=(g,max(S))$ for some $g\in G$, $H_u(U_x)$ is isomorphic to graph F). Thus, we are in the hypothesis of Theorem \ref{thm:teorema6McCord}, we get that $i$ is a weak homotopy equivalence.

\begin{figure}[h]
  \centering
    \includegraphics[scale=0.78]{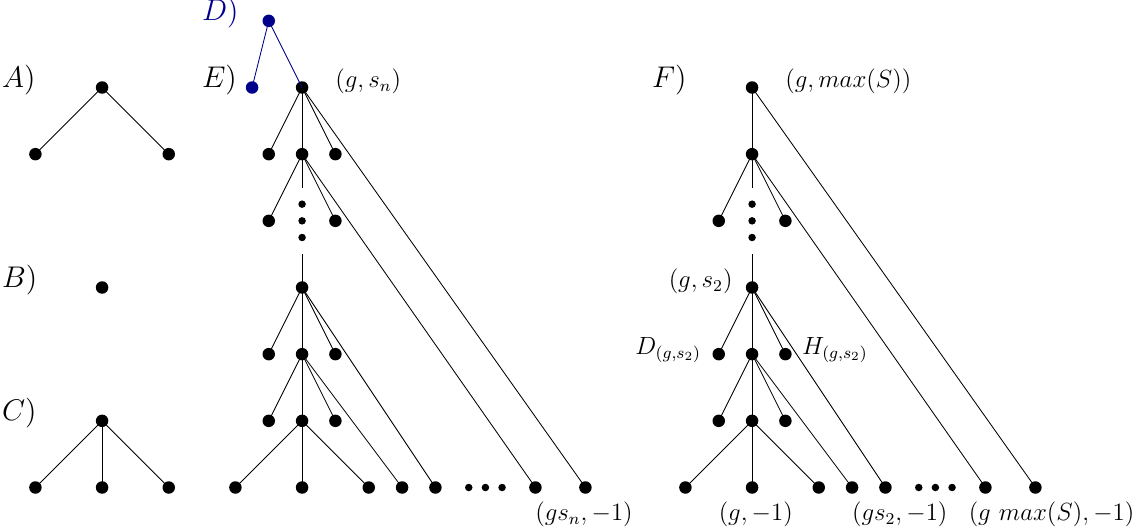}
     \caption{Graph isomorphic to $H_u(U_x)$ for different cases of $x\in\overline{X}_G$.}
     \label{figure:grafosU}
\end{figure}


We prove that $H_u(\overline{X}_G)$ is a connected $CW$-complex. We have $(g,-1)\prec (g,0)\prec (g,s_{1})\prec ... \prec (g,s_{n})\prec...$ for every $g\in G$. Repeating the same arguments used in the last part of the proof of Theorem \ref{thm:Aut(X_G)isoG}, it can be deduced the result. Therefore,  $H_u(\overline{X}_G)$ is a one-dimensional connected CW-complex, so it is homotopy equivalent to a wedge sum of circles. It remains to determine the number of circles. 

Suppose $G$ is a non-finite set. By Proposition 1A.1 \cite{hatcher2000algebraic} or Theorem 8.4.7 \cite{diek2008algebraic}, $H_u(\overline{X}_G)$ contains a maximal tree. We denote by $W_g$ the subcomplex of $H_u(\overline{X}_G)$ given by $(g,-1)\prec (g,0)\prec...\prec(g,s_{n})\prec...$. We consider a maximal tree $T$ in  $H_u(\overline{X}_G)$ containing the subcomplex $\bigcup_{g\in G}W_g$. It is clear that there is one edge of $S_{(g,\beta)}$ and another one of $T_{(g,\beta)}$ not contained in $T$ for every $(g,\beta)$, where $g\in G$ and $max(S)>\beta \geq 0$. Therefore, we have that the wedge sum of circles is at least of $2|S'||G|$ circles. But, using the arithmetic of infinite cardinals, we have that  $2|S'||G|=|G|$. On the other hand, there are $|G||S'|$ edges in $H_u(\overline{X}_G)$ not contained in $(\bigcup_{g\in G, 0\leq \beta <max(S)}H_u(S_{(g,\beta )}\cup T_{(g,\beta)}))\cup(\bigcup_{g\in G} W_g) $. Then, there are at most $2|G||S'|+|G||S'|$ edges not contained in $T$. Thus, $H_u(\overline{X}_G)$ is homotopy equivalent to $\bigvee_{\mathbb{N}}S^1$.

Suppose $G$ is finite, we know that $H_u(\overline{X}_G)$ has the homotopy type of a wedge sum of a finite number of circles. Therefore, we only need to use the Euler characteristic to determine the number of circles. The number of vertices is $v=n(r+2)+10nr$, where the first term corresponds to the number of vertices of $H_u(X_G)\subset H_u(\overline{X}_G)$ and the second one to the vertices of $H_u((S_{(g,\beta)}\cup T_{(g,\beta)})\setminus \{(g,\beta) \})$, where $g\in G$ and $0\leq \beta <max(S)$. The number of edges is $e=n(r+1)+nr+12nr$, where the two first terms correspond to the number of edges of $H_u(X_G)$ and the third term corresponds to the edges of $\bigcup_{g\in G, 0\leq \beta <max(S)} H_u(S_{(g,\beta)}\cup T_{(g,\beta)})$. Therefore, we get that the Euler characteristic  of $H_u(X_G)$ is $n-3nr$, which means that $H_u(X_G)$ is homotopy equivalent to the wedge sum of $3nr-n+1$ copies of $S^1$.


\end{proof}
\end{prop}

\begin{remark} Given a countable group $G$, the proof of Proposition \ref{prop:homotopyTypeX_G} can be adapted to $\overline{X}_G^n$ for every $n\in \mathbb{N}$. Then, it can be deduced that for every countable group $G$, $\overline{X}_G^n$ is weak homotopy equivalent to $\overline{X}_G$ for every $n\in \mathbb{N}$.
\end{remark}


\begin{prop}\label{prop:embeddingGinX_G} Given a group $G$, the McCord functor induces a natural monomorphism of groups $K:\mathcal{E}(\overline{X}_G)\rightarrow \mathcal{E}(|\mathcal{K}(\overline{X}_G)|)$.
\begin{proof}
If $[f]\in \mathcal{E}(\overline{X}_G)$, we get by Remark \ref{rem:AutEqSonIguales} that $[f]=f\in Aut(\overline{X}_G)$. Then $\mathcal{K}(f)\in Aut(|\mathcal{K}(\overline{X}_G)|)$, where $\mathcal{K}(f)$ denotes the induced map between the geometric realization of $\mathcal{K}(\overline{X}_G)$, see \cite{mccord1966singular}. It is also clear that $\mathcal{K}(f)$ defines a homotopy class such that $[\mathcal{K}(f)]  \in \mathcal{E}(|\mathcal{K}(\overline{X}_G)|)$. It is trivial to check that $K:\mathcal{E}(\overline{X}_G)\rightarrow \mathcal{E}(|\mathcal{K}(\overline{X}_G)|)$ given by $K(f)=[\mathcal{K}(f)]$ is a well-defined homomorphism of groups. We prove the injectivity of $K$, if $f,g\in \mathcal{E}(\overline{X}_G)$ with $f\neq g$, there exists $(h,\beta)\in \overline{X}_G$ with $f(h,\beta)\neq g(h,\beta)$. Therefore, $f(S_{(h,\beta)})=S_{f(h,\beta)}\neq S_{g(h,\beta)}= g(S_{(h,\beta)})$ and $\mathcal{K}(f)(|\mathcal{K} (S_{(h,\beta)})|)\neq \mathcal{K}(g)(|\mathcal{K} (S_{(h,\beta)})|)$. On the other hand, $|\mathcal{K}(S_{(h,\beta)})|$ is homotopy equivalent to $S^1$. Then, $\mathcal{K}(f)$ and $\mathcal{K}(g)$ send the same copy of $S^1$ to different copies of $S^1$ in $|\mathcal{K}(\overline{X}_G)|$. Using Proposition \ref{prop:homotopyTypeX_G}, it can be deduced that $\mathcal{K}(f)$ is not homotopic to $\mathcal{K}(g)$, so $K(f)\neq K(g)$.
\end{proof}
\end{prop}

\begin{remark} It is not difficult to check that the monomorphism of  groups defined in the proof of Proposition \ref{prop:embeddingGinX_G} is not an isomorphism of groups. For instance, we only need to consider the continuous function that interchanges $H_u(S_{(g,\beta)})$ with $H_u(T_{(g,\beta)})$ in $H_u(\overline{X}_G)$ for some $(g,\beta)$ with $g\in G$ and $0\leq \beta < max(S)$, i.e., the symmetry through $(g,\beta)$ of $H_u(S_{(g,\beta)})$ and $H_u(T_{(g,\beta)})$ in $H_u(\overline{X}_G)$ that fixes the rest of the points.
\end{remark}

For a general Alexandroff space $A$, the homomorphism of groups $K:\mathcal{E}(A)\rightarrow \mathcal{E}(|\mathcal{K}(A)|)$ given in Proposition \ref{prop:embeddingGinX_G} is not necessarily a monomorphism of groups. 
\begin{example}\label{example:notmonomorphism} Let us consider the Alexandroff space $A^c$ considered in Example \ref{example:autnoequiv}. The McCord complex $\mathcal{K}(A^c)$ of $A^c$ is a triangulation of $S^1$. Then, $\mathcal{E}(|\mathcal{K}(A^c)|)\simeq \mathcal{E}(S^1)\simeq \mathbb{Z}_2$, while $Aut(A^c)\simeq \mathcal{E}(A^c)\simeq \mathbb{Z}_2\oplus \mathbb{Z}_2$.
\end{example}

\begin{remark} In general, the image of the monomorphism of Proposition \ref{prop:embeddingGinX_G} is not a normal subgroup of $\mathcal{E}(|\mathcal{K}(\overline{X}_G)|)$. We consider the cyclic group of three elements $\mathbb{Z}_3$. By Proposition \ref{prop:homotopyTypeX_G}, $|\mathcal{K}( \overline{X}_{\mathbb{Z}_3})|$ is homotopy equivalent to $\bigvee^7_{i=1}S_i^1$, so $\mathcal{E}(|\mathcal{K}(\overline{X}_{\mathbb{Z}_3})|)\simeq  \mathcal{E}(\bigvee^7_{i=1}S_i^1)$. We take $f\in  \mathcal{E}(\overline{X}_{\mathbb{Z}_3})\simeq \mathbb{Z}_3$ of order $3$. We take $\rho \in \mathcal{E}(\bigvee^7_{i=1}S^1)$ given by $\rho(S^1_i)=S^1_{i+1}$ and $\rho(S^1_7)=S^1_{1}$, where $i=1,...,6$. It can be deduced that $K(f)$ viewed as an element of $\mathcal{E}(\bigvee^7_{i=1}S^1)$ satisfies that $K(f)(S_1^1)=S_3^1$, $K(f)(S_3^1)=S_5^1$, $K(f)(S_5^1)=S_1^1$, $K(f)(S_2^1)=S_4^{1}$, $K(f)(S_4^1)=S_6^{1}$, $K(f)(S_6^2)=S_2^{1}$ and $K(f)(S_7^1)=S_7^{1}$. Then, it is easy to check that $\rho K(f)\rho^{-1}\notin L$, where $L$ denotes $K(\mathcal{E}(\overline{X}_{\mathbb{Z}_3}))$ viewed as a subgroup of $\mathcal{E}(\bigvee^7_{i=1}S_i^1)$.

\begin{figure}[h]
  \centering
    \includegraphics[scale=0.7]{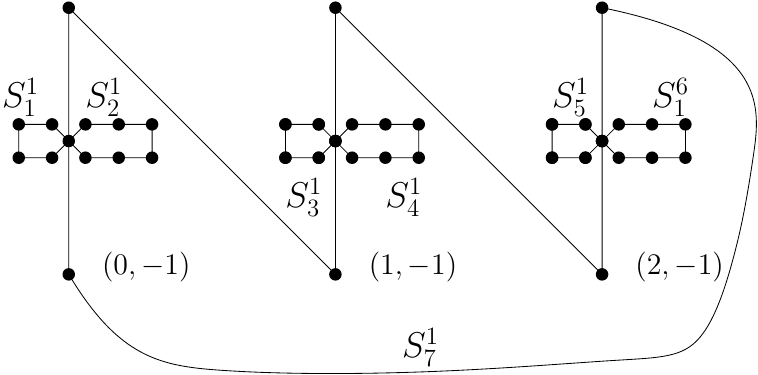}
     \caption{$H_u(\overline{X}_{\mathbb{Z}_3})$. }
     \label{figure:undirected}
\end{figure}

\end{remark}

\begin{example}\label{ex:infiniteposet} Let $\mathbb{Z}$ be the group of integer numbers with the addition, we consider $S'=\{1\}$. The Hasse diagram of $\overline{X}_\mathbb{Z}$ can be seen in Figure \ref{figure:ejemploEnterosHasse}. From the Hasse diagram, it can be deduced that $\mathbb{Z}\simeq Aut(\overline{X}_\mathbb{Z})\simeq \mathcal{E}(\overline{X}_\mathbb{Z})$ is generated by the translation to the right and to the left of the columns of the Hasse diagram, that is to say, the translations of  $\{(i,-1),(i,0),(i,1) \}$ for $i\in \mathbb{Z}$. 

\begin{figure}[h]
  \centering
    \includegraphics[scale=0.6]{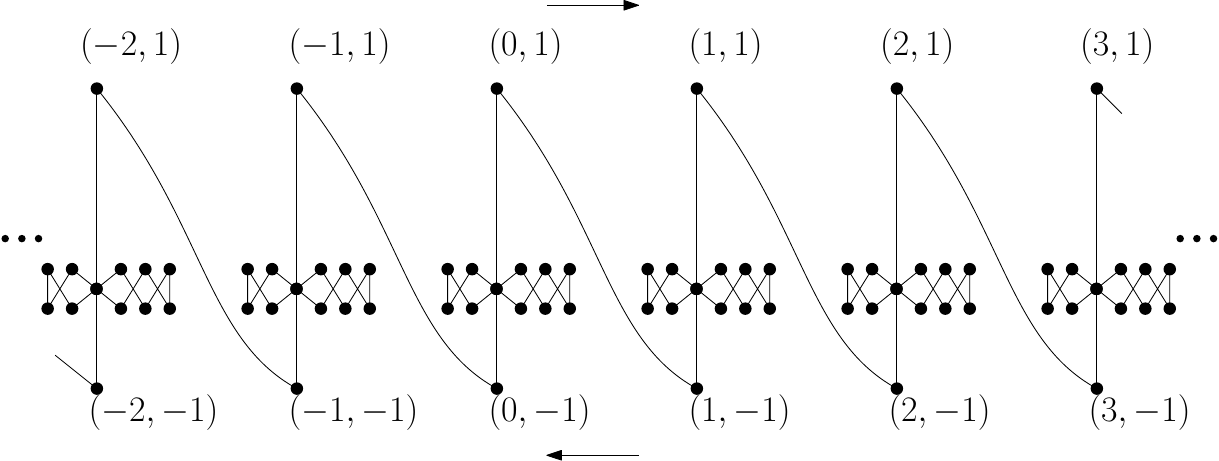}
     \caption{Hasse diagram of $\overline{X}_\mathbb{Z}$.}
     \label{figure:ejemploEnterosHasse}
\end{figure}

\end{example}

\textbf{Acknowledgments.} We would like to thank the referee for carefully reading our manuscript and for giving such valuable comments which substantially improved some previous versions of the paper.

\bibliography{bibliografia}
\bibliographystyle{acm}

\newcommand{\Addresses}{{
  \bigskip
  \footnotesize

  \textsc{ P.J. Chocano, Departamento de \'Algebra, Geometr\'ia y Topolog\'ia, Universidad Complutense de Madrid, Plaza de Ciencias 3, 28040 Madrid, Spain}\par\nopagebreak
  \textit{E-mail address}:\texttt{pedrocho@ucm.es}

  \medskip

\textsc{ M. A. Mor\'on,  Departamento de \'Algebra, Geometr\'ia y Topolog\'ia, Universidad Complutense de Madrid and Instituto de
Matematica Interdisciplinar, Plaza de Ciencias 3, 28040 Madrid, Spain}\par\nopagebreak
  \textit{E-mail address}: \texttt{ma\_moron@mat.ucm.es}

  \medskip

\textsc{ F. R. Ruiz del Portal,  Departamento de \'Algebra, Geometr\'ia y Topolog\'ia, Universidad Complutense de Madrid and Instituto de
Matematica Interdisciplinar
, Plaza de Ciencias 3, 28040 Madrid, Spain}\par\nopagebreak
  \textit{E-mail address}: \texttt{R\_Portal@mat.ucm.es}

}}

\Addresses

\end{document}